\documentclass{article}
\usepackage[utf8]{inputenc}
\usepackage{bussproofs}
\usepackage{mathrsfs}
\usepackage{amsmath}
\usepackage{amsthm}
\usepackage{amssymb}
\usepackage[dvipsnames]{xcolor}
\usepackage{stmaryrd}
\usepackage{calc}
\usepackage{cmll}
\usepackage{amsmath}
\usepackage{tikz}
\usepackage{tikz-cd}
\usepackage{url}
\usepackage{longtable}
\usetikzlibrary{arrows.meta, calc, positioning}
\usetikzlibrary{positioning,arrows,calc}
\usetikzlibrary{shapes,backgrounds}
\usepackage{stmaryrd}
\usepackage{{xargs,ifthen,xstring}}
\usepackage{authblk}
\usepackage{txfonts}
\usepackage{array}

\usepackage{tikz}
\tikzstyle{mybox} = [draw=black, very thick, rectangle, inner ysep=5pt, inner xsep=5pt]
\usepackage{caption}
\theoremstyle{definition}
\newtheorem*{fact}{Fact}

\title{Sundholm's explanation of meaning:\\logical atavism and the nature of proofs}

\author{Antonio Piccolomini d'Aragona\\

CFvW Center, University of T\"{u}bingen, T\"{u}bingen, Germany\\

\texttt{antonio.piccolomini-daragona@uni-tuebingen.de}}
\date{}

\begin{document}

\maketitle

\begin{abstract}
    I provide an overview of some of Sundholm's remarks on the history and philosophy of logic. In particular, I focus on Sundholm's proposal to explain meaning with no object-language/metalanguage distinction, and to provide a consequently contentual approach to formalisms for proofs. When applied to Gentzen's Natural Deduction in its two 1935 and 1936 variants, this triggers a reading of each variant as pointing to one of the poles in the distinction between proof-objects and proof-acts, also introduced by Sundholm. I suggest that the basis of this picture is given by a Martin-L\"{o}fian reading of the notion of analytic assertion.
\end{abstract}


\section{Introduction}

G\"{o}ran Sundholm's contributions to constructivist logic and the foundations of mathematics range from historical insights into the development of the discipline, through subtle philosophical observations and distinctions, to substantial advancements in technical issues---see \cite{klevintroduction, primierorahman}. I aim at providing a unified picture of some of the aspects of this variety of topics, which may be called \emph{Sundholm's explanation of meaning}, but which I shall refer to as \emph{Sundholm's semantics} for short.\footnote{As will be clear from what I shall be saying in this paper, the word ‘‘semantics" should not be understood here in the sense in which it is usually understood by today logicians, namely, from a mostly post-Tarskian viewpoint. In fact, it is my impression that, in Sundholm's case, the expression ‘‘semantics" is best understood in the broader Saussurean tradition \cite{saussure}, or in the sense of Kant's (J\"{a}sche's) \emph{Logik} \cite{kantlogik}.}

Sundholm's overall picture stems from two intertwined tenets at least. First, that the explanation of meaning should not be carried out under a distinction between object-language and metalanguage. Logical formalisms should be meaningful languages, and their meaning-explanation should be such as to render the axioms and rules of inference evident. This standpoint has been sometimes qualified by Sundholm as an atavistic one. Second, that logical atavism may be also applied to proof-formalisms, so that the latter, even when originally devised in a meta-theoretical perspective, be liable to a contentual reading (Section 2).

A prior role is played here by Gentzen's Natural Deduction for intuitionistic logic, which can be formalised in two ways: the familiar tree-style which Gentzen sketched in his 1935 \emph{Untersuchungen \"{u}ber das logische Schlie\ss en} \cite{gentzenuntersuchungen}, and the sequential version of it from Gentzen's 1936 \emph{Die Widerspruchsfreiheit der reinen Zahlentheorie} \cite{gentzen1936}. Although these two formalisms are commonly said to be nothing but stylistic variants of each other, it follows from Sundholm's first tenet that they are not so: once the object-language/metalanguage distinction is no longer given pride of place, formal expressions should not be understood as meaningless strings of signs for which the semantics should \emph{provide} meaning, but as expressive languages whose intended interpretation is \emph{unfolded} by the semantics. From Sundholm's second tenet it then follows that unfolding such an intended meaning is a task which the semantics can carry out, not just for the language, but also for the proof-apparatus over this language (Section 3).

So, what is the semantic unfolding of the 1935 and 1936 Natural Deduction derivations? The answer lies in Sundholm's distinction between proof-objects and proof-acts \cite{sundholm98}, which stems in turn from Sundholm's historical-philosophical research on BHK-semantics \cite{sundholm83}. Based on a Martin-L\"{o}fian reading of the analytic rendering of synthetic assertions that given propositions are true (Section 5), 1935 derivations can be seen as ‘‘gestures" towards proof-objects, while 1936 derivations are proof-traces, i.e., reifications of proof-acts. Moreover, a semantic unfolding of the purely syntactic correspondence between 1935 and 1936 derivations becomes thereby available (Sections 6 and 7).

This paper does not aim at accounting for \emph{every} aspect of Sundholm's semantics. I focus on a number of selected sources between 1983 and 2024 which I take to be more relevant for my purposes, without following a chronological order when quoting from them. I am aware that this risks concealing the \emph{evolution} of Sundholm's thought, but I leave to others or to future works the task of providing a unified and all-encompassing picture of Sundholm's semantics.

\section{Metalanguage and Bolzano's reduction}

In the logic of today, one usually proceeds as follows. First, a formal language $\mathscr{L}$ is defined as a set of strings via inductive formation-rules out of a well-rounded alphabet. Formulas of $\mathscr{L}$ are considered as meaningless objects, and it is after shifting from a purely syntactical to a genuinely semantic viewpoint that $\mathscr{L}$ comes to be endowed with meaning. Since the interpretation of the logical vocabulary is fixed throughout variations of non-logical meaning, the evaluation of a formula $A$ can be defined by induction on the logical complexity of $A$, so that $A$ becomes true or false depending on the chosen model $\mathfrak{M}$ of the non-logical terminology it involves. If truth obtains on every $\mathfrak{M}$ (possibly depending on the truth under $\mathfrak{M}$ of the elements in a set of formulas $\Gamma$), $A$ can be said to be \emph{logically} true (or a \emph{logical} consequence of $\Gamma$) since truth only stems from the logical structure of $A$ (and of the elements of $\Gamma$).

This \emph{modus operandi} is to be found e.g. in Hilbert. The idea that one and the same system could be interpreted in different ‘‘models" was a crucial move in Hilbert's \emph{Foundations of Geometry} \cite{hilbertgeometrie} for proving the relative independence and the consistency of given groups of axioms---although one should speak of models in Hilbert's case \emph{cum grano salis}, see \cite{schiemereder, klevdedekindhilber, schiemergiovannini}---while the idea that systems could or should be looked upon as meaningless strings, to be operated on mathematically, anticipated in \cite{hilbert1904}, became then functional to Hilbert's programme of achieving a finitist proof of consistency for (classical) analysis---although it is unclear whether this idea of Hilbert was more than functional to the aim of treating actual mathematics in a mathematically manageable way, not to mention the actually \emph{Inhaltliche} content of metamathematics, see e.g. \cite{detlefsen, moriconiorigini, moriconiteoriadimostrazione, sieg, smorinsky}.

The idea that formal languages come with an intended meaning was certainly endorsed by Frege, as well as by those walking in Frege's immediate footsteps, like Russell. In fact, it may have even been there in the very first metatheoretical investigations of, e.g., Carnap, G\"{o}del, and others---see e.g. \cite{awodeyreck, goldfarb, schiemerreck}. An exception is the algebraic tradition of Boole, Peirce, Schr\"{o}der, Skolem and so on, which was later on taken up partly by Hilbert himself. It is however with Lukasiewicz's and Tarski's joint seminars in Warsaw during the 1920s \cite{tarskilukasiewicz} and Tarski's \emph{The concept of truth in formalised languages} \cite{tarskitruth}, that the step is taken of keeping object-language and metalanguage apart, with the former being \emph{attributed} a meaning by the latter. A move which looked appropriate after G\"{o}del's incompleteness results \cite{godelincompleteness} and, earlier on, the failure of Frege's foundational programme because of Russell's paradox---for this reconstruction see e.g. \cite{sundholmatavism, sundholmcompleteness}

Be that as it may, it is a basic tenet of Sundholm that the \emph{metalinguistic turn} is not something we should subscribe to. In defending what he calls his \emph{atavistic position}

\begin{quote}
    that we as logicians would do well return not to our immediate ancestors, but to pre-1930 days, when interpreted formal systems still held sway, \cite[156]{sundholmatavism}
\end{quote}
Sundholm recalls what, before the advent of G\"{o}del's incompleteness theorems, the task of logicians was, i.e., not to prove

\begin{quote}
    metamathematical theorems about formal languages, metamathematically construed, but to provide a foundation for mathematics in the following way: (1) designing a sizeable formal language with an axiomatic deductive apparatus; (2) providing careful meaning-explanations for its basic (or \emph{primitive}) notions; (3) making the axioms and rules of inference evident from the meaning-explanations in question, such that (4) the resulting framework is adequate (at least in principle) for real and complex analysis. \cite[152]{sundholmatavism}
\end{quote}

In Sundholm's opinion, the insights allowed by this atavistic position are numerous. Here are the first three.

\begin{quote}
First, the object language is no longer useless, but its ‘‘expressions" will also express. [...] Second, a number of phenomena that tend to get ignored in metamathematical treatments are easily taken care of. A prime example is that of assertoric force (treated by means of the \emph{Urteilsstrich}, rather than the theorem predicate) and other pragmatic notions. [...] Third [...] in systems of Gentzen's Natural Deduction, both formation-rules and derivation rules pertain to [formulas]. Thus, they serve as formalistic simulacra of propositions [...] they also play the role of demonstrated (‘‘asserted") theorems. \cite[156]{sundholmatavism}
\end{quote}
Below, I shall extensively deal again with the second and the third advantages. As for the first, let me just remark that it is often described by Sundholm also as a distinction between understanding a formal language as a tool for \emph{speaking with}, and understanding it as an object to be \emph{spoken about}. This formulation can be found in \cite{sundholmimplicit}, but it is from another source that I would like to quote here:

\begin{quote}
  modern formal languages [...] are construed meta-mathematically [... and] can be only \emph{talked about}; they are objects of study \emph{only}, but are not intended for use [...] the ‘‘well-formed formulae" lack meaning, and as such not express. They are mathematical objects on par with other mathematical objects. \cite[552-553]{sundholmepistemicassumptions}
\end{quote}
This has a number of consequences, most notably that epistemological issues are relegated to a significantly lesser position in the logical analysis, above all when it comes to investigating the deductive behaviour of systems, and the intertwined fact that Frege's pragmatic assertion-sign $\vdash$ for propositional contents becomes a (meta)predicate indicating the derivability of formulas in a system.\footnote{In the light of Brouwer's criticism of the law of excluded middle, another consequence is the following: ‘‘after G\"{o}del's work, attempts to resuscitate Fregean logicism [...] no longer seemed viable and were abandoned: retaining classical logic as well as impredicativity, while insisting on explicit meaning-explanations that render axioms and rules of inference self-evident, simply seems to be asking too much. Thus we may jettison either meaning for the full formal language, while retaining classical logic and impredicativity, which is the option chosen by Hilbert's formalism [or] we may jettison classical logic and Platonist impredicativity, but then offer meaning explanations for constructivist language after the now familiar fashion of Heyting" \cite[554]{sundholmepistemicassumptions}. I cannot deal with this deep insight in greater detail, thus I limit myself to remarking that it would be highly interesting to explore where constructivist approaches which rely upon an object-language/metalanguage distinction, e.g. Prawitz's \cite{prawitz1973, prawitz2015}, would belong to in Sundholm's ‘‘grid". My tentative suggestion is that they jettison \emph{both} classical logic and language with content---see \cite{piccolominisundholmprawitz}. See also \cite{sundholmacenturyjudgmentinference}.}

So, we are immediately led to a fourth advantage of the atavistic view, which has to do with the fact that the metalinguistic turn is carried out under what Sundholm calls the \emph{Bolzano reduction} \cite{bolzano}. The reduction is responsible for the concealment of at least two crucial distinctions, as it splits into two steps, namely, reduction of (correctness of) assertions into (truth of) propositions, and reductions of (correctness of) inferences into (holding of) consequences. In my favourite source, Sundholm's reconstruction runs as follows. First, one observes that

\begin{quote}
    Bolzano takes the Aristotelian form of judgement and turns it into a form of content, where the contents are objectified denizens of the ideal [...] third realm [...] propositions-in-themselves serve in various logical roles, in particular as contents of mental acts and declarative sentences. \cite[9]{sundholmacenturyjudgmentinference}
\end{quote}
Based on this, one can proceed further to note that this realist move triggers another one, being the first step in the two-steps Bolzano's reduction, namely,
\begin{quote}
    the reduction of epistemological matters to the platonist \emph{an sich} notions. The first instance of this reduction concerns judgement: a judgement [...]
    \begin{center}
        proposition-in-itself $A$ is true
    \end{center}
    is correct (\emph{richtig}) if $A$ really is a truth-in-itself. \cite[10]{sundholmacenturyjudgmentinference}
\end{quote}
Finally, the second step is operated on inferences $I$ of the form
\begin{prooftree}
    \AxiomC{$J_1 \ J_2 \ ... \ J_k$}
    \UnaryInfC{$J$}
\end{prooftree}
where $J_i, J$ are judgements ($i \leq k$). Now
\begin{quote}
    the changed form of judgement transforms the inference schema $I$ into $I'$:
    \begin{prooftree}
        \AxiomC{$A_1$ is true \ $A_2$ is true \ $...$ \ $A_k$ is true}
        \UnaryInfC{$C$ is true}
    \end{prooftree}
    An inference according to $I'$ is valid if the proposition-in-itself $C$ is a logical consequence of the propositions-in-themselves $A_1, A_2, ..., A_k$. [...] In other---more modern---words, $C$ is \emph{a logical consequence} of $A_1, A_2, ..., A_k$ when the proposition-in-itself
    \begin{center}
        $(A_1 \wedge A_2 \wedge ... \wedge A_k) \rightarrow C$
    \end{center}
    is not just true, but \emph{logically} true, that is, true under all uniform variations of its non-logical parts [...] the validity of an inference is also reduced to, or, in this case perhaps better, \emph{replaced by}, something on the level of the platonist contents of judgements. \cite[11]{sundholmacenturyjudgmentinference}
\end{quote}

The price to be paid for the reduction is what, following Brentano \cite{brentanolehreurteil}, Sundholm calls \emph{blindness}. Blindness applies to both judgements and inferences. In the first case: 

\begin{quote}
    a \emph{blind} judgement, a mere guess, without any trace of justification, is a piece of \emph{knowledge} (an \emph{Erkenntnis}). \cite[11]{sundholmacenturyjudgmentinference}
\end{quote}
Whereas, in the second case,
\begin{quote}
    the inference is valid or not, irrespective of whether it does transmit knowledge from premise judgements to the conclusion judgement, solely depending on the \emph{an sich} truth-behaviour of the propositions-in-themselves that serve as contents of the judgements in questions [...] also inference can be \emph{blindly valid}, irrespective of whether it does preserve knowability from premise(s) to conclusion.\footnote{Sundholm remarks that ‘‘Bolzano also studies another consequence relation among propositions, but now restricted to the field of truths-in-themselves only, that he calls \emph{Abfolge} (grounding)" \cite[12]{sundholmacenturyjudgmentinference}. Once again, however, ‘‘the relation of grounding, \emph{which holds} in the first instance between pieces of knowledge, that is, \emph{between judgements known}, is turned into a propositional relation" \cite[12]{sundholmacenturyjudgmentinference}. Sundholm has dealt with the notion of grounding tree in \cite{sundholmgroundingtree}.} \cite[11-12]{sundholmacenturyjudgmentinference}
\end{quote}

Blindness becomes untenable when logic is understood as the science, but also the art of reasoning. Reasoning is ‘‘built up from inferences starting from (self-evident) axioms" \cite[4]{sundholmcompleteness} and, while consequence is a relation (under all variations) among truth-bearers, inferences and judgements are primarily acts:

\begin{quote}
    one \emph{draws} an inference and \emph{makes} a judgement. [...] The object [of these acts], however, is not only the objective correlate of the act. Coupled to the exercised act [...] there is also the objective signified act, that is, the trace [...]. The object (product) of an act of judgement (demonstration) is the judgement made (theorem proved). \cite[947-948]{sundholmvsrevisited}
\end{quote}
We shall come back to this. Now I would like to quote instead another piece of textual evidence concerning  why one should not be happy with epistemic blindness:

\begin{quote}
    the demonstration of the Prime Number Theorem (PNT) [...] could be formalized within NBG, the set theory of Von Neumann, Bernays and G\"{o}del. Since this theory is finitely axiomatized, we may conjoin its axioms into one proposition VNBG, and then consider the inference
    \begin{prooftree}
        \AxiomC{VNBG is true}
        \LeftLabel{(*)}
        \UnaryInfC{PNT is true}
    \end{prooftree}
    The inference (*), certainly, is truth-preserving [...]. So under the Bolzano reduction this is a valid inference, [...] but it provides no epistemic insight at all. Instead, validity of inference [...] involves preservation, or transmission, of \emph{epistemic} matters from premises to conclusion [...]. In order to validate the inference one makes the assumption that one knows the premise-judgements, or that they are being given as evident, and under this epistemic assumption\footnote{The expression ‘‘epistemic assumption" has a quite technical meaning in Sundholm's semantics. It is for example the core of a solution that Sundholm has proposed in \cite{sundholmschock} to a circularity problem raised by Prawitz, e.g. \cite{prawitzseeming}, and Martin-L\"{o}f, e.g. \cite{martinlofnormativeethics}, relative to the intertwinement of the notions of valid inference and proof. In this case, the notion also pertains to a distinction between a \emph{first-} and a \emph{third-person} perspective which one may uphold when setting up logical analysis, particularly an analysis of inferential validity---see also \cite{piccolominisundholmprawitz}. Epistemic assumptions also play a major role in Sundholm's (and van Atten's) reading of intuitionism \cite{sundholmvanatten}.} one has to make clear that also the conclusion can be made evident. \cite[556]{sundholmepistemicassumptions}
\end{quote}

Once, by means of the object-language/metalanguage split, formal languages have been devoid of an intended meaning, any deductive apparatus over the formal language loses its significance too and becomes, as the language itself, just an object to be spoken about. Via Bolzano's reduction, then, completeness and soundness theorems for given calculi come to play the role of providing an ‘‘adequacy test" of the derivability relation $\vdash$ over the semantic relation $\models$.\footnote{It should be remarked that the Bolzano reduction \emph{does not imply} the metalinguistic attitude, as is shown by the fact that Bolzano himself operated in an interpreted language.} As a matter of fact, we now conceive of these results

\begin{quote}
    as showing a deep connection between ‘‘syntax and semantics": the two turnstiles $\vdash$ and $\models$ coincide. The formulation presupposes that the \emph{relata} of the relation of consequence and inference are the same. From a contentual, ‘‘atavistic" perspective this is a very doubtful assumption. \cite[2]{sundholmcompleteness}
\end{quote}
That this is the case can be realised by looking back to Frege, to whom the introduction of the notion of turnstile is due. The ‘‘semantic" turnstile
\begin{quote}
    plays no role in Frege, and his uses of the ‘‘syntactic" turnstile is radically different from the modern one: Frege's sign serves as a pragmatic assertion indicator. \cite[552]{sundholmepistemicassumptions}
\end{quote}

\noindent Or again:

\begin{quote}
    in the Frege-like contentual paradigm [...] knowledge claims are part and parcel of the \emph{use} of the system, but not in the form of propositional operators. [...] An assertion [...] contains (implicit) claims, as to knowledge and truth, with respect to the content expressed by the sentence. \cite[194]{sundholmimplicit}
\end{quote}

If now, inspired by the atavistic spirit, we decide to go back to the good old times ‘‘when interpreted formal systems still held sway" \cite[27]{sundholmatavism}, we might wonder what the role of the semantics is to be. This cannot be, as one usually conceives of it today, that of \emph{giving} the formal language a meaning, via mappings which range onto given base-structures. For, the formal language already \emph{comes} with a meaning.\footnote{In my opinion, Sundholm's reconstruction also helps understanding why, in the very beginnings of the metatheoretical turn, when the contentual and the metamathematical approaches had not yet been kept apart sharply, there was no clear distinction between the notion of (deductive or semantic) completeness and the notion of categoricity. On this point see also \cite{awodeyreck, awodeyreck2}.} So, the semantics is expected to accomplish the task of just issuing or, as I shall be saying, of \emph{unfolding} this intended meaning, i.e.,

\begin{quote}
    in the theory of meaning, we consider interpreted formalisms, that is, formal languages, with a natural deduction derivational apparatus, where the formulae have meaning, that is, express propositions. The theory of meaning must then issue a meaning theory for the language in question. However, one would expect such a meaning theory not only to provide meaning for the terms and formulae, but also for the derivations themselves. \cite[192]{sundholm98}
\end{quote}

The last quotation points to the second basic tenet in Sundholm's semantics, namely that the atavistic approach paves the way also to a contentual interpretation of (Natural Deduction) proof-formalisms. This may look unusual to a logician trained in the metamathematical paradigm. Given Bolzano reductions, ‘‘interpreting" the alphabetic signs which the derivational apparatus operates on, seems to be enough for attaining an ‘‘interpretation" of the calculus too, as what the judgements and inferences mean is now expressed just as preservation of given semantic values among ‘‘interpreted" formulas. In particular, the ‘‘interpretation" of the proof formalism will be dealt with via suitable soundness and completeness theorems (if any).

\section{Gentzen's 1935 and 1936 Natural Deduction}

In the passage quoted at the end of Section 2, Sundholm assumes the proof-apparatus over a formal language to be a Natural Deduction one. This is because, in Sundholm's reconstruction, the Bolzanian conflation of the notions of inference and consequence is something which Gentzen's Natural Deduction was overall capable to avoid. This will be discussed in the next section. Here, I want to provide a quick reminder of some features of Natural Deduction calculi.

As is well-known, Natural Deduction is a kind of calculus introduced by Gentzen in his 1935 paper \emph{Untersuchungen \"{u}ber das logische Schlie\ss en} \cite{gentzenuntersuchungen}. In this first variant, the calculus is the one we are all familiar with now: logical constants come with rules of two types, introductions and eliminations, while derivations are defined by induction on the length of iterate applications of such rules. Both rules and derivations are given as trees where nodes are labelled by formulas from the underlying language, and where arcs instantiate applied rules. Certain rules allow for the dischargement of top-formulas, called assumptions. If all the assumptions are bound, the root-formula, i.e., the conclusion, is derived categorically while, if not, it is derived under the unbound assumptions.

I indicate by $\texttt{DER}_{35}$ the set of the 1935 Natural Deduction derivations for intuitionistic propositional logic, as given by the rules of Table 1 (brackets stand for dischargement of assumptions), and by $\vdash_{35}$ the derivability relation via the given rules. $\Gamma \vdash_{35} A$ thus means that there is $\mathscr{D} \in \texttt{DER}_{35}$ from assumptions $\Gamma$ to conclusion $A$---when $\Gamma = \emptyset$, this amounts to the categorical derivability of $A$, i.e., $\vdash_{35} A$. Negation is not a primitive, but a defined sign, namely, 
\begin{center}
$\neg A \stackrel{\text{def}}{=} A \rightarrow \bot$
\end{center}
where $\bot$ is an atomic constant symbol for absurdity. In this way, negation rules become particular instances of rules for implication.

\begin{center}
\begin{tikzpicture}
\node [mybox] (box){%
    \begin{minipage}{.98\textwidth}
\begin{prooftree}
    \AxiomC{$A \ \texttt{As}$}
    \noLine
    \UnaryInfC{}
    \AxiomC{$\bot$}
    \RightLabel{Exp}
    \UnaryInfC{$A$}
    \AxiomC{$A$}
    \AxiomC{$B$}
    \RightLabel{$\wedge_I$}
    \BinaryInfC{$A \wedge B$}
    \AxiomC{$A_1 \wedge A_2$}
    \RightLabel{$\wedge_{E, i}, i = 1, 2$}
    \UnaryInfC{$A_i$}
    \AxiomC{$[A]$}
    \noLine
    \UnaryInfC{$B$}
    \RightLabel{$\rightarrow_I$}
    \UnaryInfC{$A \rightarrow B$}
    \noLine
    \QuinaryInfC{}
 \end{prooftree}
 \begin{prooftree}
     \AxiomC{$A \rightarrow B$}
     \AxiomC{$A$}
     \RightLabel{$\rightarrow_E$}
     \BinaryInfC{$B$}
     \AxiomC{$A_i$}
     \RightLabel{$\vee_{I, i}, i = 1, 2$}
     \UnaryInfC{$A_1 \vee A_2$}
     \AxiomC{$A \vee B$}
     \AxiomC{$[A]$}
     \noLine
     \UnaryInfC{$C$}
     \AxiomC{$[B]$}
     \noLine
     \UnaryInfC{$C$}
     \RightLabel{$\vee_E$}
     \TrinaryInfC{$C$}
     \noLine
     \TrinaryInfC{}
 \end{prooftree}
    \end{minipage}
};
\end{tikzpicture}
\captionof{table}{1935 rules for intuitionistic propositional logic}
\end{center}

A second version of Natural Deduction is presented by Gentzen in his 1936 paper \emph{Die Widerspruchsfreiheit der reinen Zahlentheorie} \cite{gentzen1936}. We again have inference rules, out of which inductively defined derivations are built in tree-form. The nodes, however, are no longer formulas; following the Sequent Calculus formalism also introduced by Gentzen in \cite{gentzenuntersuchungen}, they are instead \emph{sequents}
\begin{center}
    $A_1, ..., A_n \Rightarrow B$
\end{center}
where $A_i$ and $B$ are formulas ($i \leq n$)---and $\{A_1, ..., A_n\}$ is to be understood as a multi-set. Derivations can no longer involve unbound assumptions, meaning that the leaves of the derivation-trees always instantiate an \emph{identity axiom} $A \Rightarrow A$, and the conclusion is always derived categorically. Hence, we no longer have dischargement of assumptions. The role of the 1935 assumptions is in a sense now played by the \emph{context} of the sequent, i.e., the left-hand side of it, which the right-hand side, called the \emph{end-formula}, can be understood as being entailed by. So, throughout derivations we account for operations performed both on the context and on the end-formula, which is why, besides logical rules, the calculus also involves structural rules---namely, X, C and W in Table 2 below.

I indicate by $\texttt{DER}_{36}$ the set of the 1936 Natural Deduction derivations for intuitionistic propositional logic, as per the rules of Table 2, and by $\vdash_{36}$ the derivability relation via the given rules. $\vdash_{36} \Gamma \Rightarrow A$ thus means that there is $\mathscr{D} \in \texttt{DER}_{36}$ with conclusion $\Gamma \Rightarrow A$. Once again, negation is defined as implication of $\bot$, so it will be dealt with by special instances of the implication rules.

\begin{center}
\begin{tikzpicture}
\node [mybox] (box){%
    \begin{minipage}{.98\textwidth}
\begin{prooftree}
    \AxiomC{$\Gamma, A, B \Rightarrow C$}
    \RightLabel{X}
    \UnaryInfC{$\Gamma, B, A \Rightarrow C$}
    \AxiomC{$\Gamma, A, A \Rightarrow B$}
    \RightLabel{C}
    \UnaryInfC{$\Gamma, A \Rightarrow B$}
    \AxiomC{$\Gamma \Rightarrow B$}
    \RightLabel{W}
    \UnaryInfC{$\Gamma, A \Rightarrow B$}
    \AxiomC{}
    \RightLabel{$\texttt{Id}$}
    \UnaryInfC{$A \Rightarrow A$}
    \noLine
    \QuaternaryInfC{}
\end{prooftree}
    
\begin{prooftree}
\AxiomC{$\Gamma \Rightarrow \bot$}
\RightLabel{Exp}
\UnaryInfC{$\Gamma \Rightarrow A$}
\AxiomC{$\Gamma \Rightarrow A$}
\AxiomC{$\Delta \Rightarrow B$}
\RightLabel{$\wedge_I$}
\BinaryInfC{$\Gamma, \Delta \Rightarrow A \wedge B$}
\AxiomC{$\Gamma \Rightarrow A_1 \wedge A_2$}
\RightLabel{$\wedge_{E, i}, i = 1, 2$}
\UnaryInfC{$\Gamma \Rightarrow A_i$}
\noLine
\TrinaryInfC{}
\end{prooftree}

\begin{prooftree}
    
\AxiomC{$\Gamma, A \Rightarrow B$}
\RightLabel{$\rightarrow_I$}
\UnaryInfC{$\Gamma \Rightarrow A \rightarrow B$}
    \AxiomC{$\Gamma \Rightarrow A \rightarrow B$}
    \AxiomC{$\Delta \Rightarrow A$}
    \RightLabel{$\rightarrow_E$}
    \BinaryInfC{$\Gamma, \Delta \Rightarrow B$}
    \AxiomC{$\Gamma \Rightarrow A_i$}
    \RightLabel{$\vee_{I, i}, i = 1, 2$}
    \UnaryInfC{$\Gamma \Rightarrow A_1 \vee A_2$}
    \noLine
    \TrinaryInfC{}
\end{prooftree}
\begin{prooftree}
    \AxiomC{$\Gamma \Rightarrow A \vee B$}
    \AxiomC{$\Delta, A \Rightarrow C$}
    \AxiomC{$\Theta, B \Rightarrow C$}
    \RightLabel{$\vee_E$}
    \TrinaryInfC{$\Gamma, \Delta, \Theta \Rightarrow C$}
\end{prooftree}
    \end{minipage}
};
\end{tikzpicture}
\captionof{table}{Gentzen's 1936 rules for intuitionistic propositional logic}
\end{center}

Now, it is often said that Gentzen's 1935 and 1936 versions of Natural Deduction are nothing but stylistic variants of each other. For we have, not only an extensional correspondence, in the sense that $\Gamma \vdash_{35} A$ if, and only if, $\vdash_{36} \Gamma \Rightarrow A$, but also and above all an effective, derivability-preserving way for turning any 1935 derivation into a 1936 derivation, and back.

\begin{fact}
    There is an effective mapping $\phi$ from $\texttt{DER}_{35}$ to $\texttt{DER}_{36}$ such that, when $\mathscr{D}$ witnesses $\Gamma \vdash_{35} A$, $\phi(\mathscr{D})$ witnesses $\vdash_{36} \Gamma \Rightarrow A$, and vice versa.
\end{fact}

\section{Proof-objects and proof-acts}

But the idea that the 1935 and the 1936 Natural Deduction versions are stylistic variants of each other can be upheld only when we understand those formalisms as mere ‘‘proof-games", i.e., under the metatheoretical perspective that systems are just \emph{proved about}, rather than \emph{proved in}. The \textbf{Fact} above is one such example of metatheoretical proving about.

However, if we aim to stick to the contentual approach of Sundholm's semantics, we should rather reason as follows: the 1935 and the 1936 versions amount to different formalisms, hence they \emph{express} different things, hence their semantic unfolding should be different. And since they have the same language, what must differ is the unfolding of the semantic aspects of derivations. Sundholm's semantics permits us to substantiate precisely this conclusion. Let us see how.

In \cite{sundholmvsrevisited}, Sundholm provides a thorough analysis of the distinction between \emph{implication}, \emph{conditional}, \emph{consequence} and \emph{valid inference}. Quoting from a more recent source, the distinction can be presented thus:

\begin{quote}
    (1) the \emph{implication} proposition $A \rightarrow B$, which may be \textbf{true} (or even \emph{logically true}, ‘‘in all variations"); (2) the \emph{conditional} [if $A$ is true then $B$ is true] or [...] $B$ is true \textbf{oncondition} that $A$ is true [...]; (3) the \emph{consequence} [$A \Rightarrow B$] may \textbf{hold}; (4) the \emph{inference} [$A$ is true. Therefore $B$ is true] may be valid. [...] The judgements (1)-(3) have different meaning-explanations [...] and accordingly do not mean the same, are not synonimous, while (4) indicates acts of passage. \cite[555]{sundholmepistemicassumptions}
 \end{quote}
\noindent Now, Sundholm notes first of all that when we
\begin{quote}
    consider how one would establish that (1) to (4) obtained, we see that for (1)-(3) ordinary natural deduction derivations are involved in one way or another\footnote{In particular---and in type-theoretic formalism---we need a function $b \in \texttt{Proof}(B)[x \in \texttt{Proof}(A)]$. Being a hypothetical proof of $B$ under $A$, this is already a proof-object for the conditional, whereas a proof-object for $A \rightarrow B$ will be the course-of-value of $b$, i.e., $\lambda (A, B, [x]b)$. A verification-object for $A \Rightarrow B$ is finally a higher-order function $[x]b \in (\texttt{Proof}(A))\texttt{Proof}(B)$. On this topic, besides \cite{sundholmepistemicassumptions}, from which these examples are drawn, see also \cite{klevasperus} and \cite{martinlofleiden}.}, \cite[555]{sundholmepistemicassumptions}
\end{quote}
but also, crucially, that
\begin{quote}
    with Gentzen's 1936 sequential formulation of Natural Deduction, where the derivable objects are sequents, that is consequences, [...] we get a system that can cope both with inference and consequence, \cite[555]{sundholmepistemicassumptions}
\end{quote}
precisely \emph{contra} Bolzano's reduction. At the same time, though,
\begin{quote}
it must be remembered that Gentzen wrote \emph{after} the ‘‘metalogical turn". His Natural Deduction systems [...] were designed primarily for contributions to the Hilbert programme. The formal systems are \emph{metatheoretical} in character.\footnote{For this point see also the third advantage of adopting the ‘‘atavistic" perspective pinpointed by Sundholm in \cite{sundholmatavism}, and which I quoted in Section 2. It may be also useful to recall Gentzen's position on Hilbert's programme: ‘‘a foremost characteristic of Hilbert's point of view seems to me to be the endeavour to withdraw the problem of the foundations of mathematics from \emph{philosophy} and to tackle it as far as in any way possible with methods proper to mathematics" \cite[237]{gentzenrecent}. On Hilbert's programme see, e.g., \cite{detlefsen, moriconiorigini, moriconiteoriadimostrazione, sieg, smorinsky}.} \cite[625]{sundholmvalues}
\end{quote}
Therefore, the atavistic semantics of proof-formalisms requires that metamathematical calculi like Gentzen's, and contentual systems like Frege's, be dealt with on par, so that we could either
\begin{quote}
    divest Frege systems of their content and treat them as if they were metamathematical, or [...] supply meaning explanations for the key notions in Gentzen's systems, so as to endow its object ‘‘language" with content. \cite[33]{sundholmacenturyjudgmentinference}
\end{quote}
The atavistic standpoint itself forces the second option. A step forward in this direction is taken if we turn to the distinction---already touched upon in a quotation from Section 2 above---between proof-objects, proof-acts and proof-traces.

The latter is due to Sundholm himself, who first put it forward in the context of his historical and philosophical reconstruction of the intuitionistic semantics developed by Heyting and Kolmogorov, following Brouwer's ideas---whence BHK semantics \cite{troelstravandalen}. Contrarily to the Tarskian-based framework where the meaning of the logical constants is explained by clauses which fix \emph{truth}-conditions for sentences having those constants as main sign, in BHK semantics meaning is given by \emph{proof}-conditions instead. So, for example, the meaning of $\wedge$ is given by saying that a proof of $A_1 \wedge A_2$ is a pair $\langle \pi_1, \pi_2 \rangle$ with $\pi_i$ proof of $A_i \ (i = 1, 2)$. A proof of $A_1 \vee A_2$ can be also understood as a pair, but this time of the form $\langle \pi, i \rangle \ (i = 1, 2)$, with $\pi$ proof of $A_i$, and $i$ indication of the fact that $\pi$ is a proof of the $i$-th disjunct. A proof of $A \rightarrow B$ is an (effective) function $f(x)$ such that, whenever $\pi$ is a proof of $A$, $f(\pi)$ is a proof of $B$.\footnote{Sometimes, the proof of an implication is understood as the $\lambda$-abstraction of a function as above, namely, $\lambda (A, B, [x]f)$. The BHK-clauses clearly mirror the introduction rules of Gentzen's 1935 Natural Deduction. Based on them, we can introduce non-primitive functions which mirror the elimination rules instead. Thus, left or right projection on a pair, defined by $\texttt{p}_i \ \langle \pi_1, \pi_2 \rangle = \pi_i$, can be understood as accounting for the elimination rules for conjunction. Functional application, defined by $\texttt{app} \ (\lambda (A, B, [x]f), \pi) = f(\pi)$ can be understood as elimination of implication---when $\pi$ is a proof of the antecedent. Disjunction elimination amounts to a function which selects one of two given functions from proofs of either disjunct to proofs of the conclusion, and plugs into it the left-hand element of the proof of the disjunction, i.e., $\texttt{c} x y \ (\langle \pi, i \rangle, f_1(x), f_2(x)) = f_i(\pi)$. The explosion rule, i.e., from $\bot$ to any formula $A$, is accounted for by the empty function $\emptyset_A(x)$---which is accordingly defined by the empty set of equations---since there is no proof of $\bot$. These correspondences play a crucial role in the Curry-Howard isomorphism, which I will quickly refer to in Section 6 below, see \cite{howard}.}

This proof-based approach is consonant with Heyting's claim that propositions express intentions of constructions satisfying certain conditions, while assertions signify realisations of intentions expressed by asserted propositions, namely, that intended constructions have been found \cite{heyting}---Kolmogorov would have rather spoken of problems, and of solutions of these problems respectively, see \cite{kolmogorov}. This implies that the notion of proof is in Heyting akin to the notion of construction, and Sundholm \cite{sundholm83} has pointed out that Heyting's constructions can be understood in three deeply different senses, i.e.,

\begin{itemize}
    \item[(1)] construction-acts, or
    \item[(2)] construction-objects produced by given construction-acts, or finally
    \item[(3)] reifications of construction-acts.
\end{itemize}
Against potential conflations of (1)-(3), Sundholm has argued that proofs used in clauses which fix \emph{propositional} meaning---like those for $\wedge, \vee$ and $\rightarrow$ exemplified above---should be taken as construction-objects via (2). Construction-acts as in (1) are demonstrations of \emph{assertions} that given propositions are true. These acts can be then ‘‘crystallised" as per (3), i.e., as the kind of ‘‘recipes" for performing suitable construction-acts that we normally find in mathematical textbooks or on boards in mathematics departments.

Sundholm has often come back to this tripartition, not just with the aim of interpreting BHK semantics or Heyting's views, but in the framework of a wider philosophy of proofs, or of proof-based explanations of meaning. As late as 1998 \cite{sundholm98}---but see also \cite{sundholmmanuscrito}---construction-objects, construction-acts, and reifications of construction-acts, have become, respectively, \emph{proof-objects}, \emph{proof-acts}, and \emph{proof-traces}---where the latter notion has a sort of ‘‘formal" variant referred to by Sundholm as \emph{blueprint} of a proof-act. Both proof-objects and proof-traces are objects, but they must not be confused with the \emph{object of a proof-act}, which is a special kind of assertion, i.e., a theorem proved. The full form of the latter reads ‘‘the construction(-object) $c$ is a proof(-object) of the proposition $A$" \cite[192]{sundholm98}. Such an explicit form is not proper to proved theorems only, as assertions \emph{in general} cope with the same pattern. Before addressing this issue, let me conclude by remarking that, in Sundholm's semantics, another major difference occurs between proof-objects and proof-traces, i.e., that they are \emph{not} objects of the same sort:

\begin{quote}
of the three notions of proof---act, trace, and object---the former two carry an epistemological import, whereas the last one does not: a proof-object is a mathematical object like any other, say, a function in a Banach space, or a complex contour-integral, whence, from an epistemological point of view, it is no more forcing than such objects. \cite[194]{sundholm98}
\end{quote}

\section{From synthetic to analytic assertions (and back)}

In order to understand how the semantic values of the 1935 and 1936 Natural Deduction derivations are issued in Sundholm's semantics, we must deal with Sundholm's reading of the notion of assertion. Although this can be found in various papers, I prefer to start here from the one I have been focusing on in the last part of the previous section, where we are told that
\begin{quote}
    one does not assert a proposition, but an enunciation that a proposition is true, which in fully explicit form has the shape
    \begin{itemize}
        \item[(**)] $c$ is a proof of $A$
    \end{itemize}
    where $A$ is the proposition in question. [...] The proof-object $c$ [...] is not, however, on its own a demonstration of the claim (**). \cite[200]{sundholm98}
\end{quote}
Thus, an assertion is not to be understood as the utterance of a proposition in assertive mood, but as an enunciation that a given proposition $A$ is true.\footnote{This is an important feature of Sundholm's semantics, but I cannot deal with it in detail here. The feature becomes crucial when comparing Sundholm's semantics with other akin approaches, like Prawitz's proof-theoretic semantics---which instead seems to conceive of assertions precisely as utterances of propositions in assertive mood, see e.g. \cite{prawitz2015}. The interested reader may refer to \cite{piccolominisundholmprawitz}. Concerning (a criticism to) the idea that one asserts propositions, see also \cite{martinlofleiden}.} The general form of an assertion can be thus written shortly as
\begin{center}
    $A \ \texttt{true}$.
\end{center}
However, this is only the surface appearance of an assertion. What Sundholm calls the ‘‘fully explicit shape" of it, reads instead
\begin{center}
    $c$ is a proof-object of $A$,
\end{center}
where $c$ is a proof-object of the kind of those mentioned in the BHK-clauses of Section 4. But the proof-object is not also a demonstration of the assertion. The latter is instead a proof-act, which can be reified as a proof-trace---or as the blueprint of the proof-act.

The passage from the surface appearance to the fully explicit shape of the assertion takes place under the overall constructivist spirit of Sundholm's semantics. In general, one starts with a \emph{truth-maker} analysis according to which, if proposition $A$ is true, then there is something which renders it true---in turn a refinement of the correspondence theory of truth, see e.g. \cite{dummettwtmII}. Constructively, such a truth-maker is a proof, understood as an object. So, that $A$ is true means that the class of the proof-objects of $A$ is inhabited. In Sundholm's words,
\begin{quote}
    the intuitionistic scheme
    \begin{center}
        $A$ is true = there exists a proof of $A$
    \end{center}
    can profitably be construed as a special instance of the general truth-maker analysis of truth:
    \begin{center}
        $A$ is true = there exists a truth-maker for $A$. \cite[117]{sundholmtruthmaking}
    \end{center}
\end{quote}

Sundholm remarks that the existential claim in the truth-maker analysis is not to be understood as an existential quantifier, i.e., as being of the form

\begin{center}
$(\exists x \in \mathscr{D}) \ \texttt{M}(A, x)$,
\end{center}
where $\mathscr{D}$ is a suitable domain of truth-makers, and $\texttt{M}$ a suitable (metalinguistic) truth-making relation. The reason is that

\begin{quote}
    if the existence involved here were itself to be expressed by means of an existential quantifier, firstly, in order that the quantifier be applicable, the relation of truth-making between the truth-maker and the true proposition would have to be a propositional function. Secondly, an infinite regress would arise:
    \begin{center}
        $(\exists y \in \mathscr{D}) \ \texttt{M} ((\exists x \in \mathscr{D}) \ \texttt{M}(A, x), y)$,
    \end{center}
    since the first existential quantifier $\exists y$ has yet again to be expressed in terms of the truth-condition for an existentially quantified proposition, and so on.\footnote{The notation looks slightly different in Sundholm's original writing. I have modified it for consistency with the notation I have used so far.} \cite[118-119]{sundholmtruthmaking}
\end{quote}
Precisely the same applies to the predicate $\texttt{true}$ in the surface rendering of the assertion. It cannot be understood as a metalinguistic truth-predicate, i.e., it
\begin{quote}
    cannot be propositional in nature. In other words, it is not a propositional function formulated in the language from which the sentences are taken whose meaning is explained using the predicate in question. \cite[118]{sundholmtruthmaking}
\end{quote}
The latter thesis is of course strictly intertwined with the object/act/trace structure. For, if $\texttt{true}$ in the surface rendering of assertions were propositional in nature, a proof of an assertion should end up being a proof-object of type
\begin{center}
    $\texttt{T}(A)$
\end{center}
where $\texttt{T}$ is some truth-predicate in Tarskian-style. This cannot be, however, since proof-objects are abstract mathematical objects like any others, and as such they do not carry any epistemological import. An assertion, on the contrary, conveys implicit references to knowledge, and is therefore epistemologically significant. So, it will be justified by a proof-act, possibly reified by a proof-trace, namely, by proof-notions which properly take epistemic matters into account.\footnote{The kind of regressive or circular explanations which Sundholm refers to in the passages I have just quoted constitute a major problem in Prawitz's semantics of valid arguments \cite{prawitz1973} or theory of grounds \cite{prawitz2015}. Prawitz (implictly?) endorses an object-language/metalanguage distinction, so for him both the existential in the truth-maker claims and the truth-predicate in the surface rendering of assertions are to be conceived of as propositional in nature. Likewise, an assertion is just the utterance of a proposition in assertive mood, \emph{not to be confused with the explicit ascription of truth to a propositional content}---see \cite{prawitz2015}. This confusion would for Prawitz turn an assertion into a mere \emph{metalinguistic proposition}. I have dealt with this in \cite{piccolominisundholmprawitz}.}

We can now go back from the intuitionistic understanding of the truth-maker analysis to the fully explicit form of assertions which we begun with, i.e., from
\begin{center}
    there exists a proof of $A$
\end{center}
to
\begin{center}
    $c$ is a proof-object of $A$
\end{center}
for some suitable $c$. In practically all of his writings---including the one we have been focusing on thus far in this section, i.e., \cite{sundholmtruthmaking}, but also especially in, e.g., \cite{sundholmvestiges}, \cite{sundholminferenceconsequenceimplication} and \cite{sundholmantirealism}---Sundholm refers at this point to Martin-L\"{o}f's type theory, which for him
\begin{quote}
    constitute[s] [...] the only sustained effort towards a realisation of a constructivist theory of meaning for a sizeable interpreted language serving the needs of pure mathematics on a scale comparable to that of Frege's \emph{Grundgesetze}. \cite[137]{sundholmvestiges}
\end{quote}
Thus, we know that the existential claim involved in the truth-maker analysis cannot be understood propositionally, hence it has to be explained otherwise. Following Martin-L\"{o}f, as Sundholm does, we can hence proceed as follows:
\begin{quote}
    when $A$ is a category, a general concept, then $A \ \texttt{exists}$ is a judgement. Which judgement it is, is explained by telling what knowledge is expressed by, or in, the judgement in question. [...] the judgement is explained by telling what knowledge one has to have in order to have the right to make it. Here, in order to have the right to make the judgement $A \ \texttt{exists}$, it is required that one has made some judgement $c \colon A$\footnote{This is the Martin-L\"{o}fian rendering of ‘‘$c$ is an object of type $A$".}, that is, in order to ascribe existence to a general concept one has to know an object falling under the concept in question.\footnote{This way of proceeding is not at all proper to judgements/assertions of the form $A \ \texttt{true}$ only, but to all kinds of judgements/assertions in Martin-L\"{o}f's type theory. Dealing with this, however, would lead me too far away from the purposes of this paper. The interested reader may refer to \cite{klevintroductiontypetheory} and \cite{martin-loef}. Let me also remark that, in \cite{sundholmtruthmaking}, the passage from the intuitionistic rendering of the truth-maker analysis to the fully explicit form of assertions is additionally read through the lens of the treatment of existential statements via Weyl's notion of \emph{judgement abstract}, hence of Hilbert's and Bernays' notions of \emph{incomplete communication} and \emph{partial judgement}.} \cite[123]{sundholmtruthmaking}
\end{quote}
Overall, the chain leading from the fully explicit form of assertions, through the truth-maker analysis of it, to the surface rendering in terms of truth \emph{simpliciter}, can be understood as the passage from \emph{analytic} to \emph{synthetic} assertions, where the analytic-synthetic dichotomy is here a Martin-L\"{o}fian rendering of Kant's well-known distinction \cite{kantreiner}:

\begin{quote}
   if a judgement of [the form $a \colon A$] is evident at all, then it is evident solely by virtue of the meanings of the terms that occur in it [...] an existential judgement is synthetic in Kant's terminology. [...] its evidence rests on a construction: you see, we arrive at an existential judgement, $A \ \texttt{exists}$, say, through the construction of an object which falls under the concept $A$ [...] we clearly have to go beyond what is contained in the judgement itself, namely, to the thing that exists, in order to make an existential judgement evident, and hence it must be synthetic.\footnote{In the text, Martin-L\"{o}f also refers to judgements of the form $a = b \colon A$ which, in type theory, are one of the two basic forms of judgements for elements in a type. This is because, when specifying what a type is, we must say how its (canonical) elements are formed, but also provide an identity criterion for elements in the type. As I said, dealing with this would lead me too far away, so I again limit myself to suggesting the interested reader to refer to \cite{klevintroductiontypetheory} and \cite{martin-loef}. Let me also stress that a partly different view (albeit still Martin-L\"{o}fian in spirit) on these matters is to be found in \cite{bentzenanalyticity}.} \cite[93-94]{martinloefanalytic}
\end{quote}

\section{Semantic values of derivations}

Let us now go back to the semantic values of 1935 and 1936 Natural Deduction derivations. The latter should obviously stand for proofs, but we have \emph{three} notions of proof---proof-object, proof-act, and proof-trace---thus we must choose one for each formalism. However, derivations are trees, and trees are (formal) objects, hence we may reasonably require the range of possible choices to boil down to the only two ‘‘objectual" notions at play, i.e., proof-object and proof-trace. In a nutshell, Sundholm's semantics leads to the following result: 1935 derivations point to proof-objects of propositions, whereas 1936 derivations are (formal simulacra, or blueprints of) proof-traces of proof-acts of assertions. Let us see how this can be achieved.

Based on the discussion of the previous section, assertions in Sundholm's semantics can be taken to have a surface appearance where a given proposition $A$ is ascribed with truth, i.e.,

\begin{center}
    $A \ \texttt{true}$.
\end{center}
A properly constructive unfolding of the truth-maker analysis of this yields the analytic version in fully explicit form, i.e.,
\begin{center}
    $c$ is a proof-object of $A$
\end{center}
for some suitable $c$. But this is only a special case, namely, that of \emph{categorical} assertion. Next to it, we have another kind of speech act, i.e., the \emph{conditional} truth of a proposition $A$ depending on the truth of other propositions $B_1, ..., B_n$. We can write this as follows:
\begin{center}
    $A \ \texttt{true} \ (B_1 \ \texttt{true}, ..., B_n \ \texttt{true})$.
\end{center}
One may want to understand categorical assertions as conditionals with $n = 0$. However, this would violate the order of conceptual priority in the explanation of the kind of speech act that a conditional is. As said in a quotation from Section 5, to know what an assertion is, we must say what knowledge one must have for being justified to make that assertion. Likewise, to know what a conditional is, we must say under what conditions one can licitly make that conditional claim. For doing this, we must assume that we know what kind of assertion a categorical assertion is, as the explanation must provide a criterion for the condition of the assertion $A \ \texttt{true}$ to be fulfilled upon fulfilment of the conditions of the assertions $B_i$ $\texttt{true}$ for every $i \leq n$. 

Now, since the explanation of a categorical assertion is to the effect that one must exhibit a proof-object of the proposition whose truth is being asserted, a conditional claim is to be understood as constructively justified whenever one can exhibit a means for transforming proof-objects of the $B_i$-s into proof-objects of $A$. So, the constructive unfolding of the above conditional yields a reading of it in terms of \emph{dependent proof-objects}, i.e.,
\begin{center}
    $c$ is a proof-object of $A$ provided $x_i$ is a proof-object of $B_i$ for all $i \leq n$. 
\end{center}
In order to make fully explicit that $c$'s being a proof-object of $A$ depends upon $x_i$'s being a proof-object of $B_i$ for every $i \leq n$, we could also require the dependent proof-object itself to involve the free variables $x_1, ..., x_n$, and write
\begin{center}
    $c(x_1, ..., x_n) \in \texttt{Proof}(A) \ (x_1 \in \texttt{Proof}(B_1), ..., x_n \in \texttt{Proof}(B_n))$.
\end{center}
This means that, whenever
\begin{center}
    $a_i \in \texttt{Proof}(B_i)$
\end{center}
for every $i \leq n$, then
\begin{center}
    $c(a_1, ..., a_n) \in \texttt{Proof}(A)$.
\end{center}
As Sundholm observes \cite{sundholm98}, this is tantamount to stipulating the validity of the rule
\begin{prooftree}
    \AxiomC{$a_1 \in \texttt{Proof}(B_1), ..., a_n \in \texttt{Proof}(B_n)$}
    \UnaryInfC{$c[a_1/x_1, ..., a_n/x_n] \in \texttt{Proof}(A)$}
\end{prooftree}

It is now easy to realise that a 1935 Natural Deduction derivation which witnesses $\Gamma \vdash_{35} A$ can be understood as pointing to a (dependent) proof-object of the conditional
\begin{center}
    $A \ \texttt{true} \ (B_1 \ \texttt{true}, ..., B_n \ \texttt{true})$
\end{center}
where $\Gamma = \{B_1, ..., B_n\}$. Or, in Sundholm's own words,
\begin{quote}
the derivation-trees [...] stand for certain objects, namely the (dependent) proof-objects that serve to make the proposition expressed by the end-formula of the tree a (dependently) true proposition.\footnote{This is similar to what, in \cite{martinlofleiden}, Martin-L\"{o}f calls the operation of \emph{contentualization}.} \cite[196]{sundholm98}
\end{quote}
One way for being persuaded of this is the following. Take $\mathscr{D} \in \texttt{DER}_{35}$ for $\Gamma \vdash_{35} A$, with $\Gamma = \{B_1, ..., B_n\}$, which we may indicate as
\begin{prooftree}
    \AxiomC{$B_1$}
    \AxiomC{$\dots$}
    \AxiomC{$B_n$}
    \noLine
    \TrinaryInfC{$\mathscr{D}$}
    \noLine
    \UnaryInfC{$A$}
\end{prooftree}
Now, this must stand for a proof and, if we replace its unbound assumptions with closed derivations $\mathscr{D}_i \in \texttt{DER}_{35}$ for $\vdash_{35} B_i$ for every $i \leq n$, which must also stand for proofs, we obtain a closed 1935 derivation for $\vdash_{35} A$, i.e.,
\begin{prooftree}
    \AxiomC{$\mathscr{D}_1$}
    \noLine
    \UnaryInfC{$B_1$}
    \AxiomC{$\dots$}
    \AxiomC{$\mathscr{D}_n$}
    \noLine
    \UnaryInfC{$B_n$}
    \noLine
    \TrinaryInfC{$\mathscr{D}$}
    \noLine
    \UnaryInfC{$A$}
\end{prooftree}
which must stand for a proof of $A$ too. Generalising further, the claim that 1935 derivations can be understood as (dependent) proof-objects is nothing but the effect of the Curry-Howard correspondence, whose basic insight is that the 1935 Natural Deduction derivations can be mapped onto terms of an extended typed $\lambda$-calculus, and vice versa, so as to preserve certain relations of reducibility of both derivations and terms to privileged forms, via suitable normalisation theorems---see \cite{girardproofsandtypes, hindleylercherseldin, prawitz1965}.\footnote{See footnote 8 above, which also makes clear that the mapping is based on the association of the introduction and elimination rules of 1935 Natural Deduction to, respectively, primitive and non-primitive (i.e., defined by an equation) functions for forming typed terms, and vice versa. Terms of typed $\lambda$-calculus can be roughly understood as BHK proofs of a special kind.} Whence, the semantic unfolding of the 1935 Natural Deduction formalism is to the effect that, under the Curry-Howard correspondence, a $\mathscr{D} \in \texttt{DER}_{35}$ which witnesses $\Gamma \vdash_{35} A$ with $\Gamma = \{B_1, ..., B_n\}$ may be taken to stand for a $c(x_1, ..., x_n)$ such that
\begin{center}
    $c(x_1, ..., x_n) \in \texttt{Proof}(A) \ (x_1 \in \texttt{Proof}(B_1), ..., x_n \in \texttt{Proof}(B_n))$.
\end{center}

Observe that, thus far, we needed not bother about the twofold role that formulas play in a Gentzenian formalism given by a language plus a calculus over that language. As said in a quotation of Section 2, the atavistic perspective should instead allow for a sharp distinction between formulas as expressing propositions, and formulas as standing for ‘‘assertions, that is, interpretations of derived end-formulae of proof-trees" \cite[196]{sundholm98}. The reason is that, strictly speaking, a 1935 derivation does not prove \emph{that} the proposition expressed by its end-formula is true whenever the propositions expressed by its assumption-formulas are so. Its role is that of providing a formal simulacrum of a truth-maker of a (conditional) assertion. Formulas expressing demonstrated assertions in 1935 derivations do not express in the sense that the assertions at issue are demonstrated by a suitable proof-interpretation of the derivations where they occur. Rather, they express in the sense that this proof-interpretation provides a constructive unfolding of the demonstrated assertions by exhibiting proof-objects of the asserted propositions for a fully explicit rendering of the assertions themselves. When unfolding the semantic value of a 1935 derivation, we are just showing the link between the demonstrated assertions and the propositions whose truth is being asserted, via a constructive interpretation of truth-makers as proof-objects.\footnote{Incidentally, this may also be the reason why, in approaches like Prawitz's \cite{prawitz1973, prawitz2015}, where assertions are understood as utterances of propositions in assertive mood, and where at the same time we have a constructive explanation of the notion of truth as existence of a proof, 1935 Natural Deduction derivations are indeed taken also as demonstrations of (conditional) assertions---see \cite{piccolominisundholmprawitz} for details.}

\emph{Prima facie}, the distinction between formulas as expressing propositions and formulas as expressing demonstrated assertions should be even less crucial in connection with 1936 Natural Deduction derivations since, as said, the nodes are now not formulas but sequents, and
\begin{quote}
    a sequent expresses that the consequent proposition is a consequence of the propositions that serve as antecedents. The sequent $S$:
    \begin{center}
        $A_1, ..., A_n \Rightarrow B$
    \end{center}
    expresses the relationship
    \begin{center}
        that if $A_1$ is true, then if $A_2$ is true, ..., then if $A_n$ is true, then $B$ is true
    \end{center}
    [...] when matters are as the sequent states we say that $S$ holds. \cite[197]{sundholm98}
\end{quote}
The (judgement about the) holding of a sequent is, as said in a quotation of Section 4, different from the (judgement about the) conditional truth of a given proposition under truth of other propositions. The verification-object of a sequent
\begin{center}
    $A_1, ..., A_n \Rightarrow B$
\end{center}
is not a dependent proof-object of type $B$ with variables of types $A_1, ..., A_n$, but a \emph{function} from proof-objects for $A_1, ..., A_n$ to proof-objects for $B$, i.e.,
\begin{center}
    $f \in (\texttt{Proof}(A_1))...(\texttt{Proof}(A_n))\texttt{Proof}(B)$\footnote{Let me also stress that the difference between dependent proof-objects $c(x_1, ..., x_n) \in \texttt{Proof}(B) \ (x_1 \in \texttt{Proof}(A_1), ..., x_n \in \texttt{Proof}(A_n)$ and functions $f \in(\texttt{Proof}(A_1))...(\texttt{Proof}(A_n))\texttt{Proof}(B)$ is described by Sundholm, in \cite{sundholmvsrevisited}, as the difference between functions in the Euler-Frege sense, and functions in the Riemann-Dedekind-Church sense. The distinction is discussed also by Martin-L\"{o}f as a difference between functions in the \emph{old-fashioned} sense, and functions in the \emph{modern-sense} \cite{martinlofleiden}, and used by him in \cite{martinloeftwovalues} for providing different readings of Natural Deduction which are perfectly in line with---and inspired by---Sundholm's approach. In \cite{sundholmimplicit}, Sundholm observes that ‘‘in order to form a closed consequence we need an open consequence:
    \begin{prooftree}
        \AxiomC{$c \in \texttt{Proof}(B) \ (x \in \texttt{Proof}(A))$}
        \UnaryInfC{$[x]c \in (\texttt{Proof}(A))\texttt{Proof}(B)$}
    \end{prooftree}
    is the rule in question" \cite[207]{sundholmimplicit}---where I have adapted the notation to the one I have used so far, see also footnote 6 above. In \cite{sundholmimplicit}, Sundholm provides an illuminating example for seeing the difference between the two types of function, relative to a formulation of \emph{modus ponens} as a generalised elimination in Schroeder-Heister's fashion \cite{schroeder-heister1984}. Other differences are that: while in $c(x_1, ..., x_n) \in \texttt{Proof}(B) \ (x_1 \in \texttt{Proof}(A_1), ..., x_n \in \texttt{Proof}(A_n))$ the order of the $A_i$-s is irrelevant, since substitution in $c(x_1, ..., x_n)$ can be performed simultaneously, in $f \in(\texttt{Proof}(A_1))...(\texttt{Proof}(A_n))\texttt{Proof}(B)$ the order will matter, since application may now require to respect a certain composition of the domain---see again \cite{sundholmimplicit}; while a function in the Euler-Frege sense can be $\lambda$-abstracted, and then computed out via a \emph{defined} application function $\texttt{app}$ as in footnote 6, application is \emph{primitive} in the case of a function in the Riemann-Dedekind-Church sense---see e.g. \cite{sundholmepistemicassumptions}. The idea that a proof-object of $A_1, ..., A_n \Rightarrow B$ is $f \in(\texttt{Proof}(A_1))...(\texttt{Proof}(A_n))\texttt{Proof}(B)$ seems to force $n \geq 1$. For, if $n = 0$, then a proof-object for $\Rightarrow B$ would be just an element $\in \texttt{Proof}(B)$, and it would be thus the same as a proof-object for the \emph{proposition} $B$. But, intuitively, there is a difference between $\Rightarrow B$ and $B$, which is seen when we reflect on the fact that a judgement to the effect that $\Rightarrow B$ holds is a judgement to the effect that ‘‘(I know that) it holds that $B$ is true", not simply ‘‘(I know that) $B$ is true"---see also \cite{sundholmvestiges, sundholmantirealism}.},
\end{center}
and there is no way to map a 1936 derivation into such a function, as happened instead---via Curry-Howard isomorphism---with 1935 derivations and (dependent) proof-objects of (conditional) assertions. But then, the fact that a 1936 derivation does not stand for a proof-object, and that it has nonetheless to stand for a proof of some kind, seems to imply that a 1936 derivation is nothing but the proof-trace, or better the blueprint, of a proof-act:

\begin{quote}
    a sequential derivation-tree provides a notation for an act-trace. The execution of the corresponding act of demonstration results in the claim that the consequence-relationship expressed by the end-sequent holds. \cite[198]{sundholm98}
\end{quote}
However, it is \emph{precisely} because of this that, despite the nodes of a 1936 derivation are \emph{not} formulas, when unfolding the meaning of a 1936 derivation we \emph{must} care about the clash between formulas as expressing propositions and formulas as expressing demonstrated assertions:
\begin{quote}
    from an epistemological point of view, the sequential form is the proper one in which to display natural deduction derivations: they are the representations of the \emph{demonstrations through which theorems become known}. \cite[198, italics mine]{sundholm98}
\end{quote}
To see this, we must turn to Sundholm's claim that the 1936 format is the \emph{proper} one for Natural Deduction. This can be achieved via an unfolding of the implicit pragmatic features of formulas in 1935 derivations as expressing demonstrated assertions. Let us see how---the (beautiful) unfolding is drawn from \cite{sundholmvalues}.

Assume that the end-formula of a 1935 Natural Deduction derivation is understood as expressing, not a proposition in the (contentual) language, but a demonstrated assertion in the (contentual) proof-apparatus. As said in Section 5, then, what the formula expresses is
\begin{center}
    $A \ \texttt{true}$
\end{center}
for some suitable proposition $A$. If we make this explicit in the diagrammatic rendering of a 1935 derivation from $\Gamma = \{B_1, ..., B_n\}$ to $A$ as above, we may want to move from
\begin{prooftree}
    \AxiomC{$B_1$}
    \AxiomC{$\dots$}
    \AxiomC{$B_n$}
    \noLine
    \TrinaryInfC{$\mathscr{D}$}
    \noLine
    \UnaryInfC{$A$}
\end{prooftree}
to
\begin{prooftree}
    \AxiomC{$B_1 \ \texttt{true}$}
    \AxiomC{$\dots$}
    \AxiomC{$B_n \ \texttt{true}$}
    \noLine
    \TrinaryInfC{$\mathscr{D}_1$}
    \noLine
    \UnaryInfC{$A \ \texttt{true}$}
\end{prooftree}
where $\mathscr{D}_1$ is now just $\mathscr{D}$ decorated with the explicit display of formulas-as-assertions. This first step, however, is inaccurate, since the truth of $B_1, ..., B_n$ is, not asserted, but assumed. We may make this explicit by using two distinct force indicators, $\vdash$ for asserted truth, and $\dashv$ for assumed truth, thereby obtaining
\begin{prooftree}
    \AxiomC{$\dashv B_1$}
    \AxiomC{$\dots$}
    \AxiomC{$\dashv B_n$}
    \noLine
    \TrinaryInfC{$\mathscr{D}_2$}
    \noLine
    \UnaryInfC{$\vdash A$}
\end{prooftree}
where $\mathscr{D}_2$ is $\mathscr{D}_1$ with $\vdash$ or $\dashv$ in place of $\texttt{true}$. This will not do either, though, because the \emph{dependent} truth of a proposition under truth of other propositions \emph{cannot} be taken as an outright assertion. Sundholm's example is the following derivation, using introduction and elimination of implication:
\begin{prooftree}
    \AxiomC{$[\dashv A]$}
    \noLine
    \UnaryInfC{$\mathscr{D}^*$}
    \noLine
    \UnaryInfC{$\texttt{?} \ B$}
    \RightLabel{$\rightarrow_I$}
    \UnaryInfC{$\vdash A \rightarrow B$}
    \AxiomC{$\mathscr{D}^{**}$}
    \noLine
    \UnaryInfC{$\vdash A$}
    \LeftLabel{$\mathscr{D}^{***} = \ $}
    \RightLabel{$\rightarrow_E$}
    \BinaryInfC{$\vdash B$}
\end{prooftree}
where $\texttt{?} \ B$ indicates that the semantic value of that occurrence $B$ in $\mathscr{D}^{***}$ (and in $\mathscr{D}^*$) is as of yet unsettled. But we know that---in the constructive framework---$\mathscr{D}^*$ can be taken as a $c(x)$ such that
\begin{center}
    $c(x) \in \texttt{Proof}(B) \ (x \in \texttt{Proof}(A)),$
\end{center}
whence $\texttt{?} \ B$ must be a \emph{conditional}
\begin{center}
    $B \ \texttt{true} \ (A \ \texttt{true})$.
\end{center}
Thus, all the dependent-nodes of a 1935 derivation are to be seen as conditionals. And since unbound assumptions in a 1935 derivation depend upon themselves, we find that our original $\mathscr{D}_1$ should be rather written
\begin{prooftree}
    \AxiomC{$B_1 \ \texttt{true} \ (B_1 \ \texttt{true})$}
    \AxiomC{$\dots$}
    \AxiomC{$B_n \ \texttt{true} \ (B_n \ \texttt{true})$}
    \noLine
    \TrinaryInfC{$\mathscr{D}_3$}
    \noLine
    \UnaryInfC{$A \ \texttt{true} \ (B_1 \ \texttt{true}, ..., B_n \ \texttt{true})$}
\end{prooftree}
where $\mathscr{D}_3$ is now just $\mathscr{D}_1$ where the dependencies of the conditionals are made explicit in right-hand parentheses. If we now replace
\begin{center}
    $C \ \texttt{true} \ (D_1 \ \texttt{true}, ..., D_m \ \texttt{true})$
\end{center}
by
\begin{center}
    $D_1, ..., D_m \Rightarrow C$
\end{center}
everywhere in $\mathscr{D}_3$, what we obtain is precisely a 1936 Natural Deduction derivation. Which means that
\begin{quote}
    when the deductively relevant features are made explicit in order to account also for the pragmatic interaction of assumption with assertion, standard Natural Deduction is but a variant of the Sequent Calculus version. \cite[631]{sundholmvalues}
\end{quote}
And since this stems from the unfolding of pragmatic features involved in a \emph{demonstration} process, we can conclude that 1936 derivations are proof-traces (or blueprints) of proof-acts---see also \cite{martinloeftwovalues}, as well as the more recent \cite{paginassertionassumption}. So, whenever $\mathscr{D} \in \texttt{DER}_{36}$ witnesses $\vdash_{36} \Gamma \Rightarrow A$, it can be seen as a
\begin{center}
   proof-trace for a demonstration that the sequent $\Gamma \Rightarrow A$ holds.
\end{center}

\section{Extraction procedures and objects-acts mapping}

One may at this point wonder how it happened that the unfolding of the implicit pragmatic features of a 1935 derivation, which in itself is a gesture towards a proof-object, could lead us to establish that what we had was a 1936 derivation, i.e., the blueprint of a proof-act. A first answer stems from the fact that the unfolding of the semantic values of 1936 Natural Deduction derivations\footnote{It is maybe worthwhile observing that, strictly speaking, if we take 1936 derivations as ‘‘blueprints for mental acts of knowledge [then] no \emph{objectual} semantic values are assigned" \cite[634]{sundholmvalues}. So, it is somewhat misleading to speak of the semantic value of a 1936 derivation. On the other hand, a blueprint \emph{is not} a mental act of knowledge, but a reification of it so, \emph{in some sense}, the expression is not \emph{entirely} misleading, which is why I decided to keep it in some place.}
\begin{quote}
    has been neutral with respect to the semantics used and the ensuing notion of proposition, presupposing only that the language is an interpreted one \cite[634]{sundholmvalues}
\end{quote}
---which can be easily seen also by observing that, for the unfolding described in Section 6, we can stay content with the surface appearance of assertions $A \ \texttt{true}$, without going to the (constructive) fully explicit shape $c \in \texttt{Proof}(A)$. It is when we endorse, at the language level, an explicitly intuitionistic or constructive semantics for propositions that 1935 derivations can---via Martin-L\"{o}f and Curry-Howard isomorphism---be understood as proof-objects:
\begin{quote}
    a proposition is thus explained in terms of how its canonical proof-objects may be formed and when two such objects are equal canonical proofs. \cite[634]{sundholmvalues}
\end{quote}
But the phenomenon has also another, to my mind deeper, explanation, which has to do with the fact that Sundholm's semantics allows for the nice property that proof-objects and proof-traces are in a sense invertible with each other. That is, we are provided here with a semantic rendering, not only of 1935 and 1936 derivations as such, but also of the (allegedly metatheoretic) \textbf{Fact} of Section 3 that an effective procedure exists for extracting a 1935 derivation out of a 1936 one, and vice versa. At the semantic level this should yield---and \emph{in fact} yields in Sundholm's semantics---that given a proof-object, we can extract from it the proof-trace of the proof-act delivering it as a truth-maker of the proposition involved in the proven assertion, and vice versa.

The correspondence is not proved \emph{out of} the \textbf{Fact} of Section 3, but rather the other way around. The effective translatability of 1935 derivations into 1936 derivations, and vice versa, is not a \emph{meta}-property of \emph{meaningless} formalisms, but the residual \emph{blueprint} of a proper issuing of the \emph{intended meaning} of the formalisms. This is what Sundholm's following question leads to:

\begin{quote}
    from a derivation in the Gentzen standard format, we easily find a corresponding sequential derivation. What significance does this metamathematical observation carry on the level of interpreted, meaningful formalisms? \cite[198]{sundholm98}
\end{quote}
The answer obtains by considering, first of all, a correct, i.e., demonstrable assertion
\begin{center}
    (*) $\ c$ is a proof-object of proposition $A$.
\end{center}
Since (*) is demonstrable, a demonstration of it must exist. And since (*) is an \emph{analytic} assertion, (*) itself must contain the relevant information for building such a demonstration up. And
\begin{quote}
    this is where the proof-object $c$ comes in handy. [...] [It] comprises information as to (i) the rule according to which it was formed, (ii) the inferentially relevant components out of which the proposition $A$ was built, and (iii) proof-objects for the premises of the rule in (i) according to which the object $c$ was formed. \cite[198-199]{sundholm98}
\end{quote}
Thus, $c$ allows for an \emph{extraction} of a proper demonstration of (*). Sundholm's example for seeing this is as follows: suppose that (*) has the form
\begin{center}
    $\langle \pi_1, \pi_2 \rangle$ is a proof-object of $A_1 \wedge A_2$.
\end{center}
By meaning explanation of $\wedge$, $\pi_i$ is a proof-object of $A_i \ (i = 1, 2)$. Thus
\begin{quote}
    the information contained in such a proof-object allows us to span---from below---a formal structure, which can be read with content as a \emph{proof-trace} $P$ with (*) as its conclusion, \cite[199]{sundholm98}
\end{quote}
namely,
\begin{prooftree}
    \AxiomC{$\pi_1$ is a proof-object of $A_1$}
    \AxiomC{$\pi_2$ is a proof-object of $A_2$}
    \BinaryInfC{$\langle \pi_1, \pi_2 \rangle$ is a proof-object of $A_1 \wedge A_2$}
\end{prooftree}
Now, decoding in turn the information contained in $\pi_1$ and $\pi_2$ respectively,
\begin{quote}
    the structure can be extended upwards by means of suitable similar steps, that depend on the particular forms of the objects and propositions involved. [...] All steps [...] will be completely deterministic, owing to the careful design of the proof-object machinery. Eventually, because of the presupposed correctness of (*), the various branches of the edifice that has been generated will break off in suitable axioms [...]. In this fashion, then, a formal structure is spanned. When considered in the direction opposite to that of its generation, that is, from top to bottom, it can be read as a blue-print for acts of demonstration with the claim (*) as its object. \cite[199]{sundholm98}
\end{quote}
Thus, given a proof-object $c$ in a correct assertion that $c$ is a proof-object of proposition $A$, Sundholm's semantics states that one can always extract from $c$ a proof-trace of a proof-act of the assertion itself.\footnote{A number of remarks are in order here. First, the reader should refer to Section 4, where the different forms of proof-objects for formulas/propositions of different logical forms are given, following BHK semantics. Second, Sundholm's example in \cite{sundholm98} does not employ the (polymorphic) notation $\langle \pi_1, \pi_2 \rangle$, but the (monomorphic) notation where types of proof-objects are made explicit, e.g., $(\langle \pi_1, \pi_2 \rangle, A_1, A_2)$---see also footnote 6 above. This is required for the extraction process to be deterministic in nature, but I preferred to leave this out not to burden the notation excessively. Third, Sundholm's reasoning seems to refer to the phenomenon of decidability of assertions $a : A$---possibly under assumptions---in Martin-L\"{o}f's type theory via \emph{type-checking algorithms}. For example, in a recent paper, and within the context of a discussion of the relation between incompleteness phenomena and Martin-L\"{o}f's type theory, Sundholm states that ‘‘the \emph{logic of analytic judgments} of the form $a : A$ is complete and decidable (by type-checking)" \cite[48]{sundholmcompleteness}---on this, see also \cite{martinloefanalytic}. The connection between analyticity and decidability via type-checking (and meaning-explanation) has been recently revised by Bentzen in \cite{bentzenanalyticity}. Bentzen's conclusion may not be compatible with Sundholm's claim that, given a proof-object in a correct assertion, we are always capable of extracting from it a proof-trace for demonstrating that assertion. However, it is my impression that this claim by Sundholm should be understood as a sort of \emph{semantic requirement}, which must be satisfied in \emph{ideal} semantic conditions. It seems to me that, in a strictly Sundholmian perspective, if a semantics is incapable of coping with the extraction-principle, then this semantics is simply ill-behaving.} That the inverse is also the case is almost trivial: for any proof-trace of a proof-act of the assertion that $c$ is a proof-object of $A$, we just pick $c$ out, thereby extracting from the proof-trace the proof-object of $A$.

\section{Conclusion}

Let us recap the \emph{fil rouge} connecting the claims which I started my overview of Sundholm's semantics with, i.e., that the meaning-explanation should be carried out under a contentual perspective where the object-language/meta-language distinction plays no role, and that such an approach allows for a contentual reading of (Natural Deduction) proof-formalism, to the final outcomes discussed in Sections 6 and 7, namely, that 1935 derivations are gestures towards proof-objects, while 1936 derivations are proof-traces, and that one can always extract a proof-object from a proof-trace, and vice versa:

\begin{enumerate}
\item the idea that one should not subscribe to the metalinguistic perspective implies two things: first, that the role of a semantics must be that of \emph{unfolding an intended meaning}, rather than that of \emph{attributing the meaning which the formalisms come without} and, second, that the 1935 and the 1936 versions of Natural Deduction must mean different things;

\item Natural Deduction, in particular in its 1936 variant, comes to play a major role, since it is overall capable to avoid the Bolzano reduction, by allowing for a distinction between conditionals, consequence, and inference. So, we must find the semantic values of the Natural Deduction derivations, and this is done in two steps:

\item first, by distinguishing proof-objects, proof-acts, and proof-traces and

\item second, by providing a constructive, Martin-L\"{o}fian rendering of assertions as existence of proof-objects as truth-makers of propositions, based on the observation that existence and truth cannot be here propositional in nature;

\item we come to the actual semantic unfolding of the 1935 and 1936 derivations, i.e., the former are gestures towards proof-objects of propositions, while the latter are proof-traces (or better blueprints) of proof-acts;

\item finally, we can thereby also provide a semantic interpretation of the purely syntactic inter-translatibility of 1935 and 1936 derivations, in terms of possibility of extracting proof-objects from proof-traces, and vice versa.
\end{enumerate}

As is seen, Sundholm's semantics as reconstructed here amounts to a bold theoretical framework whose most systematic parts are drawn from observations pertaining to the history and the philosophy of logic. The access point and path I chose may not be the only possible ones, but I believe that, whatever strategy one adopts, one will end up with roughly the same results---with respect to the topics that this paper focuses on. Also, the fact that one can tackle Sundholm's semantics from many different perspectives and following many different trails is, I think, a strong selling point of the approach.

This also means, though, that much more could be said than what I have done here. There are many issues which I could touch upon only quickly, or could not touch upon at all---for example, Sundholm's notion of \emph{epistemic assumption}, for which see footnote 3 above, or Sundholm's dialogical reading of inferences, see, e.g., \cite{sundholmschock}. Since much more work remains to be done, I would hence like to conclude by saying, with Nietzsche, ‘‘there is another world to discover---and more than one! On to the ships, you philosophers!" \cite[163]{nietzsche}.

\paragraph{Acknowledgements} I am indebted to Ansten Klev and G\"{o}ran Sundholm for stimulating discussions and precious insights on earlier drafts of this paper. This work has been financially supported by the grant PI 1965/1-1 for the DFG project \emph{Revolutions and paradigms in logic. The case of proof-theoretic semantics}.

\bibliographystyle{abbrv}
\bibliography{bibliography}

\end{document}